\documentclass{amsart}
\usepackage{epsfig,color,mathrsfs}
\title[The Perron-Frobenius Theorem for Homogeneous, Monotone Functions]{The Perron-Frobenius Theorem for Homogeneous, Monotone Functions}
\author{St\'ephane Gaubert}
\address{INRIA, Domaine de Voluceau,
B.P.~105, 78153 Le Chesnay C\'edex, France.}
\email{Stephane.Gaubert@inria.fr}
\author{Jeremy Gunawardena}
\address{Bauer Center for Genomics Research, Harvard University, 7 Divinity Avenue, Cambridge, MA 02139, USA.}
\email{jgunawardena@cgr.harvard.edu}

\def\N{{\mathbb N}}     
\def\R{{\mathbb R}}     
\newcommand{\ha}[1]{\hat{#1}}
\newcommand{\mrm}[1]{\text{\rm #1}}
\newcommand{\ra}{\rightarrow}

\newcommand{\Rp}{\R^{+}}
\newcommand{\Rim}{{\R \cup \{-\infty\}}}
\newcommand{\mychi}{\chi}
\renewcommand{\top}{\mathsf{t}}
\renewcommand{\bot}{\mathsf{b}}
\newcommand{\supnorm}[1]{{\| #1 \|_\infty}}
\newcommand{\hilbert}[1]{{\| #1 \|_{\mathsf H}}}
\newcommand{\E}{{\mathcal E}}
\newcommand{\G}{{\mathcal G}}
\newcommand{\cP}{{\mathscr{P}}}
\newcommand{\set}[2]{\{#1\mid\,#2\}}
\newcommand{\comp}{\circ}
\newcommand{\dhil}{{{\sf d}_{\sf H}}}

\newtheorem{theorem}{Theorem}
\newtheorem{cor}{Corollary}
\newtheorem{lemma}{Lemma}
\newtheorem{proposition}{Proposition}
\theoremstyle{remark}

\begin{document}
\begin{abstract}
If $A$ is a nonnegative matrix whose associated directed graph is strongly connected, the Perron-Frobenius theorem asserts that $A$ has an eigenvector in the positive cone, $(\Rp)^n$. We associate a directed graph to any homogeneous, monotone function, $f: (\Rp)^n \ra (\Rp)^n$, and show that if the graph is strongly connected then $f$ has a (nonlinear) eigenvector in $(\Rp)^n$. Several results in the literature emerge as corollaries. Our methods show that the Perron-Frobenius theorem is ``really'' about the boundedness of invariant subsets in the Hilbert projective metric. They lead to further existence results and open problems.
\end{abstract}
\maketitle

\section{Introduction and statement of main results}

This introduction provides an overview of the paper. We state the main results but defer some definitions to later sections. 

\subsection{The Perron-Frobenius theorem}
\label{s2-pft}

The classical Perron-Frobenius theorem may be stated as follows (see \cite[Chapter~2]{bp94} for more background). Let $A$ be a $n \times n$ nonnegative matrix. The graph associated to $A$, $\G(A)$, is the directed graph with vertices
$1, \cdots, n$ and an edge from $i$ to $j$ if, and only if, $A_{ij} \not= 0$. A directed graph is said to be {\em strongly connected} if there is a directed path between any two distinct vertices. The strong connectedness of $\G(A)$ is equivalent to requiring that $A$ is an irreducible matrix. Let $\Rp = \{x \in \R \;|\; x > 0 \}$ denote the positive reals. 

\begin{theorem}[Classical Perron-Frobenius theorem]
\label{theo-pf}
If $\G(A)$ is strongly connected then $A$ has an eigenvector in $(\Rp)^n$, unique up to a scalar multiple, whose associated eigenvalue is the spectral radius of $A$.
\end{theorem}

Now let $f: (\Rp)^n \ra (\Rp)^n$ be a self-map of the positive cone which satisfies the following properties. 
\begin{eqnarray}
\label{homo}
\forall \lambda\in \Rp \;\mrm{and}\; \forall x \in (\Rp)^n, \;\;f(\lambda x) =\lambda f(x) \enspace, \\
\label{mono}
\forall x,y \in (\Rp)^n, \;\;x \leq y \implies f(x) \leq f(y) \enspace .
\end{eqnarray}
The first property is called {\em homogeneity}; the second {\em monotonicity}. Here, $x \leq y$ denotes the product ordering on $\R^n$: $x \leq y$ if, and only if, $x_i \leq y_i$ for all $1 \leq i \leq n$. If $A$ is a $n \times n$ nonnegative matrix, the map $f(x)=Ax$ satisfies both properties, although only those matrices with no zero row---nondegenerate matrices---yield self-maps of the positive cone.

If $u \in \Rp$ and $J \subseteq \{1, \cdots, n\}$, let $u_J \in (\Rp)^n$ denote the vector defined by
\begin{equation}
\label{eq-uj}
(u_J)_i = \left\{\begin{array}{ll}
u & \mrm{if $i \in J$} \\
1 & \mrm{if $i \not\in J$} \enspace.
                 \end{array}\right.
\end{equation}
If $f: (\Rp)^n \ra (\Rp)^n$ is any homogeneous, monotone function, then $f_i(u_J)$ is a monotone function of $u$ for any $J$. Define the associated graph of $f$, $\G(f)$, to be the directed graph with vertices $1, \cdots, n$ and an edge from $i$ to $j$ if, and only if, 
\begin{equation}
\label{eq-ga}
\lim_{u\to\infty}f_i(u_{\{j\}}) = \infty \, . 
\end{equation}
If $f$ is a linear map, represented by the nonnegative matrix $A$, then it is clear that $\G(f)$ is identical to $\G(A)$. A vector $x \in (\Rp)^n$ is a (nonlinear) {\em eigenvector} if $f(x) = \lambda x$, for some (nonlinear) {\em eigenvalue} $\lambda \in \Rp$. We prove the following.

\begin{theorem}[Generalised Perron-Frobenius theorem]
\label{theo-new}
Let $f: (\Rp)^n \ra (\Rp)^n$ be any homogeneous, monotone function. If $\G(f)$ is strongly connected then $f$ has an eigenvector in $(\Rp)^n$.
\end{theorem}

The eigenvalue, which is the same for any eigenvector in $(\Rp)^n$---see (\ref{e-lambda}), is characterised in Proposition~\ref{prop-collatz-wielandt} below. Readers familiar with the linear theory will recognise this as a Collatz-Wielandt property. (We note that this is sufficient to obtain the spectral radius statement in Theorem~\ref{theo-pf}). The eigenvalue may be considered as a spectral {\em radius}, in a limited sense, by extending $f$ to the boundary of $(\Rp)^n$, which may always be done continuously, \cite{bns03a}, and considering eigenvectors lying in the boundary; see \cite[Theorem~3.1(1)]{nuss86}. When $f$ satisfies a suitable convexity condition, a simplified construction of $\G(f)$ is possible; see Proposition~\ref{prop-convex} in \S\ref{s2-graph}.

Consider as an example the homogeneous, monotone function
\begin{equation}
f(x) = \left(\begin{array}{c}
a\sqrt{x_1x_2} \wedge a'\sqrt{x_2x_3}\\
b\sqrt{x_2x_3} \vee b'\sqrt{x_3x_1}\\
cx_1 \vee c'x_3
\end{array}\right)\enspace,
\label{e-gex}
\end{equation}
where $a,a',b,b',c,c'$ are arbitrary parameters in $\Rp$. We use $\vee$ and $\wedge$ as infix notations for $\max$ and $\min$, respectively. Using (\ref{eq-ga}), it is easy to see that $\G(f)$ is the graph
\[
\begin{picture}(0,0)%
\includegraphics{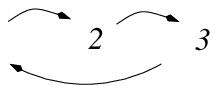}%
\end{picture}%
\setlength{\unitlength}{1973sp}%
\begingroup\makeatletter\ifx\SetFigFont\undefined%
\gdef\SetFigFont#1#2#3#4#5{%
  \reset@font\fontsize{#1}{#2pt}%
  \fontfamily{#3}\fontseries{#4}\fontshape{#5}%
  \selectfont}%
\fi\endgroup%
\begin{picture}(2268,743)(1996,-1720)
\put(1996,-1374){\makebox(0,0)[lb]{\smash{\SetFigFont{10}{12.0}{\rmdefault}{\mddefault}{\itdefault}{\color[rgb]{0,0,0}1}%
}}}
\end{picture}

\]
which is strongly connected. Hence, by Theorem~\ref{theo-new}, $f$ has an eigenvector, independently of the values of the parameters. In contrast to the linear case, the eigenvector need not be unique, even up to a positive scalar multiple. For instance, the homogeneous, monotone function 
\begin{equation}
\label{eq-xunq}
f(x_1,x_2) = (x_1 \vee x_2/2 \,,\, x_1/2 \vee x_2) \enspace,
\end{equation}
has eigenvalue $1$ and eigenspace $\{x \in (\Rp)^2 \;|\; x_1/2 \leq x_2 \leq 2x_1\}$. Here, $\G(f)$ is not just strongly connected; it is fully connected, with an edge between any pair of vertices. In the linear case, this would correspond to a positive matrix.

There is an extensive literature on homogeneous, monotone functions; see \cite{mxpnxp0,nuss88} for references. Several results somewhat similar to Theorem~\ref{theo-new} have appeared, the most relevant among them, to the best of our knowledge, being those of Nussbaum, \cite[Theorem~4.1]{nuss86}, and Amghibech and Dellacherie, \cite{ad94}. 

Nussbaum makes the following definition. Let $f: (\Rp)^n \ra (\Rp)^n$ be a homogeneous, monotone function and $A$ an $n \times n$ nonnegative matrix. $f$ is said to have $A$ as an {\em incidence matrix with respect to being power-bounded below}, \cite[Definition~4.2]{nuss86}, if, whenever $A_{ij} \not= 0$, there exists $c \in \Rp$ and a stochastic vector $\sigma \in (\Rp)^n$ with $\sigma_j > 0$, such that, for all $x \in (\Rp)^n$, 
\begin{equation}
\label{eq-nuss}
f_i(x) \geq cx_1^{\sigma_1} \cdots x_n^{\sigma_n} \enspace.
\end{equation}
Nussbaum shows that if $A$ is irreducible then $f$ has an eigenvector in $(\Rp)^n$, \cite[Theorem~4.1]{nuss86}. (The statement of \cite[Theorem~4.1]{nuss86} uses the additional hypothesis that $f$ is superadditive and then deduces that the eigenvector is unique. However, as pointed out in \cite{nuss86}, superadditivity is not required to show existence.) Since (\ref{eq-nuss}) and the condition $\sigma_j > 0$ together imply (\ref{eq-ga}), $\G(f)$ must be strongly connected whenever $A$ is irreducible. Hence, \cite[Theorem~4.1]{nuss86} is a corollary of Theorem~\ref{theo-new}. It is also not difficult to check that there is no irreducible $3 \times 3$ matrix which can be an incidence matrix for example (\ref{e-gex}). The result of Amghibech and Dellacherie follows from Theorem~\ref{theo-new} and Proposition~\ref{prop-collatz-wielandt} of the present paper but it will be convenient, for notational reasons, to defer a discussion to \S\ref{s2-graph}.

Among classical results that are corollaries of Theorem~\ref{theo-new} we should mention those of Bather, \cite[Theorem~2.4]{bat73}, in stochastic control, and Zijm, \cite[Theorem~3.4]{zijm84}, in mathematical economics. We leave the deductions to the interested reader. The ``max-plus'' version of the Perron-Frobenius theorem, \cite[Theorem~3.23]{fsch92}, is also an immediate corollary. Max-plus algebra and idempotent analysis, \cite{max95,idem,km97,ms92}, while not explicit in this paper, were crucial in stimulating the ideas behind it.

The methods needed to prove Theorem~\ref{theo-new} require several concepts---orbits, invariant subsets, super-eigenspaces, Hilbert's projective metric---which do not appear in its statement. These give rise to a new way of thinking about the eigenvector problem, to further results on existence of eigenvectors and to several open problems. We discuss these in the remainder of this Introduction.

\subsection{Invariant sets and Hilbert's projective metric}
\label{s2-hilbert}
Hilbert's {\em projective metric}, defined in (\ref{def-hilbert}), is a function  $\dhil: (\Rp)^n \times (\Rp)^n \ra \Rp \cup \{0\}$ which satisfies all the conditions of a metric except that $\dhil(y,z) = 0$ if, and only if, $y = \lambda z$ for some $\lambda \in \Rp$. It gives rise to a metric on the projective space of lines in $(\Rp)^n$, from which property it gets its name. An important observation---see (\ref{nxp-hilbert})---is that any homogeneous, monotone function $f: (\Rp)^n \ra (\Rp)^n$ is nonexpansive with respect to the Hilbert metric: for all $x,y \in (\Rp)^n$,
\begin{equation}
\label{eq-nxphil}
\dhil(f(x), f(y)) \leq \dhil(x,y) \enspace.
\end{equation}

The following characterisation is at the heart of the present paper. Up to a trivial modification, it is a special case of a theorem of Nussbaum, \cite[Theorem~4.1]{nuss88}, stated in the context of nonexpansive functions acting on cones in Banach spaces. To keep our account both elementary and self-contained, we state and prove the special case expressed in Theorem~\ref{theo-orbits}. Recall that an orbit of $f$ is any set of the form $\{ f^k(x) \;|\; k \in \N \}$ for some $x \in (\Rp)^n$.

\begin{theorem}
\label{theo-orbits}
Let $f: (\Rp)^n \to (\Rp)^n$ be a homogeneous, monotone function. $f$ has an eigenvector in $(\Rp)^n$ if, and only if, some (and hence all) orbits of $f$ are bounded in the Hilbert projective metric.
\end{theorem}

If $f$ has an eigenvector, $f(x) = \lambda x$, then the orbit of $x$ is $\{\lambda^k x\}$. This is evidently bounded in the Hilbert metric: its diameter is zero. If any orbit is bounded, then (\ref{eq-nxphil}) shows that all orbits must be so too. The force of Theorem~\ref{theo-orbits} lies in the assertion that some orbit being bounded is sufficient for the existence of an eigenvector. 

Theorem~\ref{theo-orbits} provides a simple prescription for determining the presence of an eigenvector: find an invariant subset of $(\Rp)^n$---a subset $A$ such that $f(A) \subseteq A$---which is bounded in the Hilbert projective metric. Since $A$ is invariant, it must decompose into orbits, all of which must then be bounded. For a specific function, an appropriate invariant subset may present itself naturally. However, there are invariant subsets that can be defined for all homogeneous, monotone functions.

Let $f: (\Rp)^n \ra (\Rp)^n$ be a homogeneous, monotone function and let $\lambda \in \Rp$. The {\em super-eigenspace} of $f$ corresponding to $\lambda$, $S^\lambda(f)$, is defined by $S^\lambda(f) = \{ x \in (\Rp)^n \;|\; f(x) \leq \lambda x \}$. Properties (\ref{homo}) and (\ref{mono}) show that $S^\lambda(f)$ is invariant under $f$. Note that $S^\lambda(f) \not= \emptyset$, if $\lambda$ is large enough. 

\begin{theorem}
\label{theo-main-new3}
If $f: (\Rp)^n \to (\Rp)^n$ is a homogeneous, monotone function such that $\G(f)$ is strongly connected, then all super-eigenspaces of $f$ are bounded in the Hilbert projective metric.
\end{theorem}

Super-eigenspaces are interesting, as invariant subsets, because their boundedness in the Hilbert metric can be determined by such combinatorial (graph-theoretic) constructions. These provide an effective mechanism for using Theorem~\ref{theo-orbits}; so effective, indeed, that its use can be distilled into Theorem~\ref{theo-new}, which follows immediately from Theorems~\ref{theo-orbits} and~\ref{theo-main-new3} and makes no mention of orbits, invariant subsets or Hilbert's projective metric. 

Super-eigenspaces appear at least as far back as Krein and Rutman's famous result on the existence of eigenvectors for linear functions acting on cones in Banach spaces, \cite{krut48}. This was inspired by topological fixed point theory, particularly Brouwer's Theorem and Schauder's infinite dimensional generalisation. While broadly applicable, these give little control over the location of a fixed point: the eigenvector may lie in the boundary of the cone, not in its interior. Krasnoselskii's classic text, building on Krein and Rutman's work, makes extensive use of lattice structures as an alternative to fixed point theory; see, for instance, \cite[Theorem~4.1]{kra64}. The principal difference between the present paper and that of the Russian school lies in the use of the Hilbert projective metric, for the origins of which, see \cite{nuss88}.

The ideas introduced in this sub-section suggest many further questions, some of which are discussed in sub-sections \S\ref{s2-ind} and \S\ref{s2-slice} which follow.

\subsection{Indecomposability}
\label{s2-ind}
The boundedness of all super-eigenspaces can be characterised by a combinatorial property, at the expense of comparing subsets of vertices. This is related to ideas developed in the mathematical economics literature, which we discuss further below.

We use the notation introduced in (\ref{eq-uj}). A homogeneous, monotone function $f: (\Rp)^n \ra (\Rp)^n$ is {\em decomposable} if there is a partition $I \cup J = \{1, \cdots, n\}$, $I \cap J = \emptyset$, such that, $\forall i \in I$,
\begin{equation}
\label{eq-decomp}
\lim_{u \ra \infty} f_i(u_J) < \infty \enspace.
\end{equation}
$f$ is {\em indecomposable} if it is not decomposable.

\begin{theorem}
\label{theo-main-new4bis}
All the super-eigenspaces of a homogeneous, monotone function, $f: (\Rp)^n \to (\Rp)^n$, are bounded in the Hilbert projective metric if, and only if, $f$ is indecomposable. 
\end{theorem}

It is not in fact necessary to consider all partitions of $\{1, \cdots, n\}$. There is an alternative test which relies on the recursive construction of directed graphs $\G^k(f)$. We give details in \S\ref{s2-aggregated}. For the moment, let us say that $\G^1(f)$ coincides with $\G(f)$, while, for $k \geq 2$, $\G^k(f)$ is obtained from $\G^{k-1}(f)$ by aggregating its strongly connected components. The process stabilises in the sense that there is a least integer $N \leq n$, such that $\G^N(f)$, $\G^{N+1}(f)$, $\cdots$ are isomorphic. Let $\G^\infty(f) = \G^N(f)$. (If $f$ satisfies the same convexity conditions as for Proposition~\ref{prop-convex} in \S\ref{s2-graph}, then the aggregation process stops at or before $k=2$; see Proposition~\ref{prop-convex2} in \S\ref{s2-aggregated}.)

\begin{theorem}
\label{theo-main-new4}
A homogeneous, monotone function $f: (\Rp)^n \to (\Rp)^n$ is indecomposable if, and only if, $\G^{\infty}(f)$ is strongly connected. 
\end{theorem}

Consider the following homogeneous, monotone function
\begin{equation}
f(x) = \left(\begin{array}{c}
x_1 \vee 2/(1/x_2 + 2/x_3 + 1/ x_4) \\
7x_3 \wedge x_4 \\
8\sqrt[3]{x_1x_2x_4} \\
x_3 \vee x_4 
\end{array}\right)\enspace .
\label{eq-example2}
\end{equation}
$\G(f)$ is easily seen to be
\[
\begin{picture}(0,0)%
\includegraphics{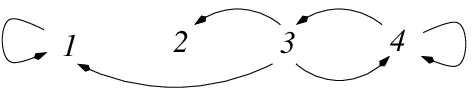}%
\end{picture}%
\setlength{\unitlength}{1973sp}%
\begingroup\makeatletter\ifx\SetFigFont\undefined%
\gdef\SetFigFont#1#2#3#4#5{%
  \reset@font\fontsize{#1}{#2pt}%
  \fontfamily{#3}\fontseries{#4}\fontshape{#5}%
  \selectfont}%
\fi\endgroup%
\begin{picture}(4474,770)(1444,-1720)
\end{picture}

\]
which is {\em not} strongly connected, so we cannot apply Theorem~\ref{theo-new}. However, it will follow from the discussion in \S\ref{s2-aggregated} that $\G^\infty(f)$, which coincides with $\G^4(f)$, is strongly connected. Hence, by Theorems~\ref{theo-orbits}, \ref{theo-main-new4bis} and~\ref{theo-main-new4}, $f$ has a positive eigenvector. Indeed, $f(u)=2u$, with $u=(1,2,8,4)^T$. 

Oshime, \cite{osh83}, following earlier work of Morishima, \cite[Appendix]{mor64}, has introduced {\em non-sectional} functions. These are close to {\em indecomposable} functions but yield a unique eigenvector in $(\Rp)^n$, \cite[Theorem~8]{osh83}. Uniqueness is important in the mathematical economics tradition; properties (\ref{homo}) and (\ref{mono}) reflect consumption behaviour in an economy and uniqueness presumably conveys a comforting sense of stability. Dietzenbacher, \cite{die94}, summarises work in this direction, going back to Solow and Samuelson's 1953 paper, \cite{ss53}. A similar concern with uniqueness can be seen in the population biology literature; see \cite{nuss86,nuss89} for references. As we observed with example (\ref{eq-xunq}), uniqueness of the eigenvector cannot be expected, even when Theorem~\ref{theo-new} can be applied. Unlike the classical linear case, conditions which guarantee uniqueness are quite restrictive. 

On the subject of uniqueness, we note that example (\ref{eq-xunq}) cannot be a contraction in the Hilbert metric or any other metric, or else the Banach Contraction Theorem would yield a unique eigenvector, \cite[Chapter~2]{gk90}. (Eigenvectors correspond to fixed points, either projectively or by rescaling the function.) This rules out the possibility, suggested by (\ref{eq-nxphil}), of using the Contraction Theorem in the manner of Birkhoff's proof of Perron's theorem on positive matrices, \cite{bir57}; see also \cite{kop82}. It can be brought into play more subtly, using the vanishing discount method of stochastic optimal control, but this leads to rather different results, which relate the existence of eigenvectors to the asymptotic dynamics of $f$, \cite{wodes,kohl80}.

\subsection{Slice spaces and recession functions}
\label{s2-slice}
Super-eigenspaces are not the only invariant subsets that can be defined for all homogeneous, monotone functions. They have dual {\em sub-eigenspaces}, $S_\mu(f) = \{ x \in (\Rp)^n \;|\; \mu x \leq f(x) \}$. The duality comes from the functional $f \ra f^-$, such that $f^-(x) = f(x^{-1})^{-1}$. (We denote by $x^{-1}$ the vector $(x_1^{-1}, \cdots, x_n^{-1})$.) Results for super-eigenspaces have corresponding dual results for sub-eigenspaces. For instance, in the dual graph to $\G(f)$, there is an edge from $i$ to $j$ if $\lim_{u \ra 0^+} f_i(u_{\{j\}}) = 0$. We generally leave it to the reader to state these dual results.

The intersection of a super-eigenspace and a sub-eigenspace forms a {\em slice space}, $S^\lambda_\mu(f) = S^\lambda(f) \cap S_\mu(f)$, which is also clearly invariant under $f$. Slice spaces are more powerful than either super-eigenspaces or sub-eigenspaces, since the boundedness of $S^\lambda_\mu(f)$ follows from that of either $S^\lambda(f)$ or $S_\mu(f)$. The following result gives a sufficient condition for the boundedness of slice spaces. 

Let $f: (\Rp)^n \ra (\Rp)^n$ be a homogeneous, monotone function. Suppose that, for each $x \in (\Rp)^n$, the limit 
\[ \hat{f}(x) = \lim_{k \ra \infty} f(x_1^k, \cdots, x_n^k)^{1/k} \]
exists, the $k$-th root being applied to each component of $f$. The function $x \ra \hat{f}(x)$ is then seen to be another homogeneous, monotone function, which we refer to as the {\em recession function} associated to $f$. (Recession functions do not always exist, as we show by example in \S\ref{s2-slice2}, but they do for all reasonable functions, including, in particular, all the other examples discussed in this paper.) Recession functions always have the unit vector as an eigenvector:
\[ \hat{f}(1, \cdots, 1) = (1, \cdots, 1) \,. \]

\begin{theorem}
\label{theo-main-new5}
Let $f: (\Rp)^n \ra (\Rp)^n$ be a homogeneous, monotone function whose associated recession function, $\hat{f}$, exists. Suppose that $\hat{f}$ has only the unit vector as an eigenvector, up to a positive scalar multiple. Then all the slice spaces of $f$ are bounded in the Hilbert projective metric.
\end{theorem}
We note that $\hat{f}$ encodes both the behaviour of $f$ at $\infty$ and at $0$. (It is a self-dual notion: if $\hat{f}$ exists, then so does $\widehat{f^-}$ and $\widehat{f^-} = \hat{f}^-$.) Here, the {\em uniqueness} of the eigenvector for $\hat{f}$ implies the {\em existence} of an eigenvector for $f$. 

Consider the homogeneous, monotone function
\begin{equation}
\label{e-ill}
f(x)=\left(\begin{array}{c}
x_2 \vee x_3\\
(x_1 \vee  x_2)\wedge x_3\\
(x_2 \vee  x_3)\wedge x_1
\end{array}\right)\enspace .
\end{equation}
It is easy to check that $f=\hat{f}$ and that if $\hat{f}(x)=x$ then $x_1=x_2=x_3$. It follows from Theorem~\ref{theo-main-new5} that all slice spaces of $f$ are bounded in the Hilbert metric. However, for all $\nu\geq 1$, $(\nu,1,1)\in S^{1}(f)$ and $(\nu,1,\nu)\in S_1(f)$, which shows that all the non-empty sub- and super-eigenspaces of $f$ are unbounded. Of course, in this example the existence of an eigenvector is trivial. However, $f$ can be altered in interesting ways without changing its recession function. Consider the homogeneous, monotone function
\begin{equation}
\label{e-ill2}
g(x)=\left(\begin{array}{c}
a_1x_2 + b_1x_3\\
((a_2x_1 +b_2x_2)^{-1} + c_2x_3^{-1})^{-1} \\
((a_3x_2 +b_3x_3)^{-1} + c_3x_1^{-1})^{-1}
\end{array}\right)\enspace ,
\end{equation}
where the $a_1, \cdots, a_3$, $b_1, \cdots, b_3$ and $c_1, \cdots, c_3$ are all in $\Rp$. It is easy to see that $\hat{g} = f$. It follows that $g$ has an eigenvector in $(\Rp)^3$. 

It remains an open problem whether the boundedness of all slice spaces can be determined by combinatorial or graph-theoretic constructions as in Theorems~\ref{theo-main-new3} and~\ref{theo-main-new4}.

\subsection{Summary and conclusions}
Functions which are homogeneous and monotone on the positive cone provide a natural generalisation of nonnegative matrices and have a correspondingly wide range of application: in population biology, \cite{nuss89}, mathematical economics, \cite{die94,mor64}, discrete event systems, \cite{fsch92,mxpnxp0}, idempotent analysis, \cite{idem,km97,ms92}, stochastic control and game theory, \cite{ag03,kol98,rs01}, nonlinear potential theory, \cite{del90}, among others. A fundamental property of such functions is their nonexpansiveness in the Hilbert projective metric, (\ref{eq-nxphil}). From the perspective of the present paper, functions which are homogeneous and monotone differ from functions which are only nonexpansive, in having invariant subsets, such as super-eigenspaces, which can be defined uniformly for all such functions. Our main contribution here has been to develop methods for determining the boundedness of these invariant subsets in the Hilbert metric---Theorems~\ref{theo-main-new3}, \ref{theo-main-new4bis}, \ref{theo-main-new4} and \ref{theo-main-new5}---from which the existence of a positive eigenvector follows. (This last implication, proved independently in Theorem~\ref{theo-orbits}, is due to Nussbaum, who demonstrates it in the context of nonexpansive functions on cones in Banach spaces, \cite[Theorem~4.1]{nuss88}.) This gives a new perspective on the classical Perron-Frobenius theorem, which is generalised in Theorem~\ref{theo-new} to any homogeneous, monotone function. We show that the Perron-Frobenius theorem is ``really'' about the boundedness of certain invariant subspaces in the Hilbert projective metric.

There are two main directions to explore in the light of the present results. Firstly, they apply to functions in which {\em all} invariant subsets of a particular type are bounded. Such functions have a quality of {\em stability}, in that graphs like $\G(f)$ depend only the behaviour of $f$ at $\infty$: the existence, or not, of an eigenvector is invariant under perturbations which preserve (\ref{homo}) and (\ref{mono}) and do not alter the divergence in (\ref{eq-ga}). (We observed this behaviour with (\ref{e-gex}) as well as (\ref{e-ill}) and (\ref{e-ill2}).) A more delicate question arises for functions in which some, but not all, non-empty, invariant subsets of a particular type are bounded. The existence of an eigenvector then depends sensitively on the structure of $f$ and on the values of its parameters. While the methods developed here for estimating boundedness---see Lemma~\ref{lem-new}---are sometimes helpful, we lack general results to handle this situation.

Secondly, all the concepts used here can be defined for functions acting on cones in Banach spaces. It is particularly interesting, in the light of Nussbaum's general theorem on orbits, to ask if the perspective and methods of the present paper can be extended to this more general setting. We know of no way to do this.

A preliminary account of some of these results appeared in \cite{xevmhf}. A special case of Proposition~\ref{prop-collatz-wielandt} below (which is the same as Proposition~1 of \cite{xevmhf}) is proved in \cite[Formula (5.2)]{bm03}. Theorem~\ref{theo-fhat} has been used in \cite{ag03}. We are grateful to the reviewer of \cite{xevmhf} for directing our attention to the mathematical economics literature; to Cormac Walsh for pointing out the significance of \cite[Theorem~4.1]{nuss88}, which we had not observed when \cite{xevmhf} was written, and to Roger Nussbaum, for comments pertaining to Theorem~\ref{theo-fhat}.

\section{The eigenvalue and the Collatz-Wielandt property}
\label{s1-evalue}

\subsection{The additive framework}

In the rest of the paper we shall work in the additive framework. We first introduce this and then explain our reasons.

The whole space, $\R^n$, can be placed in bijective correspondence with the positive cone, $(\Rp)^n$, via the mutually inverse bijections $\exp: \R^n \ra (\Rp)^n$ and $\log: (\Rp)^n \ra \R^n$, where $\exp(x) = (\exp(x_1), \cdots, \exp(x_n))$, for $x \in \R^n$, and $\log(y) = (\log(y_1), \cdots, \log(y_n))$, for $y \in (\Rp)^n$. If $f: (\Rp)^n \ra (\Rp)^n$ is any self-map of the positive cone, let $\E(f): \R^n \ra \R^n$ denote the function $\E(f)(x) = \log(f(\exp(x)))$. This induces a bijective functional between self-maps of $(\Rp)^n$ and self-maps of $\R^n$. Clearly, $\E(fg) = \E(f)\E(g)$, so that the dynamics of $f$ on $(\Rp)^n$ and $\E(f)$ on $\R^n$ are equivalent. Let $f: \R^n \ra \R^n$. The properties of $\exp$ and $\log$ show that homogeneity and monotonicity of $\E^{-1}(f)$ are equivalent, respectively, to the following properties of $f$:
\begin{eqnarray}
\label{homo2}
\forall h \in \R \;\mrm{and}\; \forall x \in \R^n, \;\;f(x + h) = f(x) + h\enspace, \\
\label{mono2}
\forall x,y \in \R^n, \;\;x \leq y \implies f(x) \leq f(y) \enspace .
\end{eqnarray}
The partial order is, as before, the product ordering on $\R^n$. We use in (\ref{homo2}) the following {\em vector-scalar} convention: if, in a binary relation or operation, a vector and a scalar appear together, the relation is taken to hold, or the operation is applied, to each component of the vector. We are accustomed to this with $\lambda x$, where $\lambda \in \R$ and $x \in \R^n$, but it is useful to extend it to, for instance, $(x + h)$ and $x \leq h$. If $x \in \R^n$ and $h \in \R$, these mean, respectively, that, for all $1 \leq i \leq n$, 
\[ (x+h)_i = x_i + h \;\;\mbox{and}\;\; x_i \leq h \enspace. \]

Gunawardena and Keane refer to functions $f: \R^n \ra \R^n$, which satisfy (\ref{homo2}) and (\ref{mono2}), as {\em topical}, \cite{ctnxp}. We shall use this terminology here and reserve the qualifiers {\em homogeneous and monotone} for functions on $(\Rp)^n$ which satisfy (\ref{homo}) and (\ref{mono}). It is equivalent to work either additively, with topical functions, or multiplicatively, with homogeneous, monotone functions. However, certain constructions are more intuitive on one side than the other: we find the Hilbert projective metric easier to work with additively, as in (\ref{def-suphil}). 

We choose, as a matter of technical convenience, to make our proofs in the additive framework. All the concepts and results introduced in the Introduction have equivalent additive formulations, for which we use the same names. For instance, the topical function $f: \R^n \ra \R^n$ has an eigenvector, $u \in \R^n$, if $f(u) = u + h$ for some $h \in \R$; the duality functional on topical functions takes $f(x)$ to $-f(-x)$; and so on. For the most part, we give new definitions and state separate versions of the theorems above, so that the rest of the paper should be self-contained. As a general rule, we leave it to the reader to formulate any dual results.

\subsection{Nonexpansiveness and cycle times}
\label{s2-nxp}

A key property of topical functions is their nonexpansiveness with respect to certain norms and metrics. Let $\top,\bot: \R^n \ra \R$ be defined as (``top'') $\top(x) = x_1 \vee \cdots \vee x_n$, and (``bottom'') $\bot(x) = -\top(-x) = x_1 \wedge \cdots \wedge x_n$. Note that both $\top$ and $\bot$ satisfy the analogues of (\ref{homo2}) and (\ref{mono2}) for functions $\R^n \ra \R^1$. The supremum norm and the Hilbert semi-norm on $\R^n$ are given, respectively, by
\begin{equation}
\label{def-suphil}
\supnorm{x} = \top(x) \vee -\bot(x) \;\;\mbox{and}\;\; \hilbert{x} = \top(x) - \bot(x) \enspace.
\end{equation}
In the notation of (\ref{eq-nxphil}), if $y,z \in (\Rp)^n$, then
\begin{equation}
\label{def-hilbert}
\dhil(y,z) = \hilbert{\log(y) - \log(z)} \enspace.
\end{equation}
The multiplicative version of the supremum metric is known as Thompson's {\em part metric} on $(\Rp)^n$, \cite[Chapter~1]{nuss88}.

An elementary application of (\ref{homo2}) and (\ref{mono2}), \cite[Proposition~1.1]{ctnxp}, shows that a function $f:\R^n\to\R^n$ is topical if, and only if,\[ \forall x,y \in\R^n, \top(f(x) - f(y)) \leq \top(x - y) \enspace. \]
(This provides some justification for the term {\em topical}.) We see immediately that a topical function is nonexpansive with respect to both the supremum norm and the Hilbert semi-norm: $\forall x,y \in \R^n$,
\begin{eqnarray}
\label{nxp-supremum}
\supnorm{f(x)-f(y)} & \leq & \supnorm{x-y} \\
\label{nxp-hilbert}  
\hilbert{f(x)-f(y)} & \leq & \hilbert{x-y} \enspace.
\end{eqnarray}
In fact, as observed by Crandall and Tartar~\cite{ct80}, if $f$ is homogeneous, then it is monotone if, and only if, it is nonexpansive in the supremum norm, \cite[Proposition~1.1]{ctnxp}.

The nonexpansiveness property (\ref{nxp-supremum}) implies that all trajectories of $f$ are asymptotically the same to within a constant:
\begin{equation} 
\label{eq-scales}
f^k(x) = f^k(y) + O(1) \;\;\mrm{as $k \ra \infty$}\enspace.
\end{equation}
(We mean by this that the function $\N \ra \R$ given by $k \ra \supnorm{f^k(x) - f^k(y)}$ is bounded as $k \ra \infty$.) Taking for $x$ an eigenvector of $f$ with associated eigenvalue $\lambda$, it follows from \eqref{homo2} and \eqref{mono2} that $f^k(x)=k\lambda +x$. Hence, using~\eqref{eq-scales},
\begin{equation}
\label{e-lambda}
\lambda = \lim_{k\to \infty } f^k(y)/k \enspace,
\end{equation}
for all $y\in \R^n$. In particular, the eigenvalue $\lambda$
is unique.

More generally, nonexpansiveness allows us to deduce that certain averages on trajectories are independent of the trajectory and give rise to functionals on the space of topical functions.
For instance, an elementary argument using (\ref{homo2}) and (\ref{mono2}) shows that the sequence $\top(f^k(0))$ is sub-additive,
\[ \top(f^{k+l}(0)) \leq \top(f^k(0)) + \top(f^l(0)) \enspace. \]
A dual inequality holds for $\bot(f^k(0))$. It follows that the sequences $\top(f^k(x)/k)$ and $\bot(f^k(x)/k)$ both converge as $k \to\infty$ and that the limits are independent of $x$. The {\em upper cycle-time} of $f$, 
$\overline{\mychi}(f) \in \R$, \cite[Definition~2.1]{ctnxp}, is defined as
\begin{equation}
\label{def-top}
\overline{\mychi}(f) = \lim_{k \ra \infty}\top(f^k(x)/k) \enspace.
\end{equation}
Dually, the {\em lower cycle-time} is $\underline{\mychi}(f) = \lim_{k \ra \infty}\bot(f^k(x)/k)$. We observe from this that, for any $k \in \N$, 
\begin{equation}
\label{eq-mychik}
\overline{\mychi}(f^k) = k\overline{\mychi}(f) \enspace,
\end{equation}
and similarly for $\underline{\mychi}$. We will make use of of this below. 

The existence of the {\em cycle-time vector} of $f$, $\mychi(f) = \lim_{k \ra \infty}f^k(x)/k \in \R^n$, is more delicate. It does not always exist, \cite[Theorem~3.1]{ctnxp}, and one of the central problems in the subject is to characterise those topical functions for which it does. 

One class of functions for which it does are functions with an eigenvector. In this case, all the functionals discussed here collapse to the eigenvalue,
since, by~\eqref{e-lambda}, 
\begin{subequations}
\begin{eqnarray}
\label{eq-topbot}
\underline{\mychi}(f) = \lambda = \overline{\mychi}(f) \;\;\mrm{and} \\
\label{eq-mychi}
\mychi(f) = (\lambda, \cdots, \lambda) \enspace.
\end{eqnarray}
\end{subequations}

\subsection{Sub-eigenspaces, super-eigenspaces and slice spaces}
Let $f: \R^n \to \R^n$ be a topical function and $\lambda,\mu \in \R$. The {\em super-eigenspace}, $S^\lambda(f)$, {\em sub-eigenspace}, $S_\mu(f)$, and {\em slice space}, $S_\mu^\lambda(f)$, are defined by
\begin{eqnarray*}
S^\lambda(f) & = & \set{x \in \R^n}{f(x) \leq \lambda +x} \\
S_\mu(f) & = & \set{x \in \R^n}{\mu +x \leq f(x)} \\
S_\mu^\lambda(f) & = & \set{x \in \R^n}{\mu + x \leq f(x) \leq \lambda +x} \enspace.
\end{eqnarray*}
It follows immediately from (\ref{homo2}) and (\ref{mono2}) that all such spaces are invariant subsets. When working with super-eigenspaces we leave it to the reader to formulate the dual results for sub-eigenspaces. It is easy to see that for any topical functions $f,g: \R^n \to \R^n$ and any $\lambda,\mu\in\R$, 
\begin{equation}
\label{eq-mono2}
\lambda\leq\mu \implies S^\lambda(f) \subseteq S^\mu(f)\enspace,
\end{equation}
Let $\Lambda(f) \subseteq \R$ denote the set of those $\lambda$ for which the corresponding super-eigenspace is non-empty: $\Lambda(f) = \set{\lambda \in \R}{S^\lambda(f)\neq \emptyset}$. It follows from (\ref{eq-mono2}) that $\Lambda(f)$ must be an interval of the form $(-\infty,\infty)$, $(a,\infty)$ or $[a,\infty)$. 

The first form can be ruled out. Suppose that $f: \R^n\to\R^n$ is a topical function and that $f(x) \leq \lambda + x$ for some $x \in \R^n$ and some $\lambda \in \R$. Using (\ref{homo2}) and (\ref{mono2}), $f^k(x) \leq k\lambda + x$. Hence, 
\[ \top(f^k(x)/k) \leq \lambda + \top(x/k) \enspace. \]
Letting $k \ra \infty$, we deduce the following lemma.
\begin{lemma}
\label{lemma-Lambda}
If $f:\R^n\to\R^n$ is a topical function then either $\Lambda(f)=(a,\infty)$ or $\Lambda(f)=[a,\infty)$, where $\overline{\mychi}(f) \leq a$.
\end{lemma}
Both possibilities can occur. It follows from (\ref{eq-topbot}) and Proposition~\ref{prop-collatz-wielandt} below that if $f$ has an eigenvector, $f(x) = \lambda + x$, then $\Lambda(f) = [\lambda,\infty)$. If $f = \E(A)$ where $A$ is the nonnegative matrix below
\[
\left(\begin{array}{cc}
      1 & 1 \\
      0 & 1
      \end{array}\right)
\]
then it is easy to see that $\Lambda(f) = (0,\infty)$.

\subsection{A key lemma and the Collatz-Wielandt property}

It remains to identify the $a$ that appears in Lemma~\ref{lemma-Lambda}. We shall show that, in fact, $a = \overline{\mychi}(f)$. This requires the following simple but crucial lemma. We extend the infix notation $\wedge$ and $\vee$ componentwise to vectors in $\R^n$.

\begin{lemma}
\label{lemma-byk}
Let $f:\R^n \to \R^n$ be a topical function and let $k$ be any positive integer. If $S^\lambda(f^k) \neq \emptyset$, then $S^{\lambda/k}(f)\neq \emptyset$.
\end{lemma}
\begin{proof}
If $S^\lambda(f^k)\neq \emptyset$, then $f^k(x) \leq \lambda +x $ for some $x\in \R^n$. Let 
\[ y = x \wedge (f(x)-\lambda/k) \wedge \cdots \wedge (f^{k-1}(x)-(k-1)\lambda/k) \enspace. \]
Using (\ref{homo}) and (\ref{mono}) we see that
\begin{eqnarray*} 
f(y) & \leq & f(x) \wedge (f^2(x)-\lambda/k) \wedge \cdots \wedge (f^k(x)-(k-1)\lambda/k) \\
     & \leq & f(x) \wedge (f^2(x)-\lambda/k) \wedge \cdots \wedge (x+\lambda/k) \\
     & =    & y + \lambda/k \enspace.
\end{eqnarray*}
Thus, $y\in S^{\lambda/k}(f)\neq\emptyset$.
\end{proof}

Lemma~\ref{lemma-byk} allows us to give the following characterisation of $\overline{\mychi}(f)$ which is an additive version of the classical Collatz-Wielandt formula, \cite[\S1.3]{min88}. 

\begin{proposition}[Generalised Collatz-Wielandt formula]
\label{prop-collatz-wielandt}
Let $f:\R^n\to\R^n$ be a topical function. Then, 
\begin{equation}
\inf \Lambda(f) = \inf_{x \in \R^n} \top (f(x)-x) = \overline{\mychi}(f) \enspace.
\label{eq-collatz}
\end{equation}
\end{proposition}
\begin{proof}
Let $a = \inf \Lambda(f)$. Since $f(x) \leq x+ \lambda$ if, and only if, $\top (f(x)-x)\leq \lambda$ the first equality in (\ref{eq-collatz}) follows easily. Lemma~\ref{lemma-Lambda} has already shown that $\overline{\mychi}(f) \leq a$. Now choose $\epsilon>0$. For sufficiently large $k$, $f^k(0) \leq (\overline{\mychi}(f) +\epsilon)k$. Hence, $S^{(\overline{\mychi}(f)+\epsilon)k}(f^k)\neq\emptyset$. By Lemma~\ref{lemma-byk}, $S^{\overline{\mychi}(f)+\epsilon}(f)\neq\emptyset$. Hence, $a \leq \overline{\mychi}(f)+\epsilon$. Since $\epsilon$ was chosen arbitrarily, $a \leq \overline{\mychi}(f)$ and so $a = \overline{\mychi}(f)$. 
\end{proof}

\subsection{Coordinates which realise cycle times} (The results of this sub-section are not required in the rest of the paper). As a byproduct of the method of Lemma~\ref{lemma-byk}, we can answer in the affirmative a conjecture of Gunawardena and Keane, \cite[Conjecture~2.1]{ctnxp}. This shows that although the cycle time vector may not exist in general, some coordinate must converge to $\overline{\mychi}$. (By duality, a similar assertion holds for $\underline{\mychi}$.) We begin with a more precise statement; the conjecture is Corollary~\ref{cor-gk} below.

\begin{theorem}
\label{theo-dyn}
Let $f: \R^n \ra \R^n$ be a topical function, and let $x\in \R^n$. There exists $1 \leq i \leq n$, such that, for all $k \in \N$,
\begin{eqnarray} 
x_i + k \overline{\mychi}(f) & \leq & f_i^k(x) \label{e-new}
\enspace . 
\end{eqnarray}
\end{theorem}

\begin{proof} 
Let $g= f - \overline{\mychi}(f)$. 
By (\ref{homo2}), $\overline{\mychi}(g) = 0$. For all $k \in \N$ let
$y(k) = x \wedge g(x) \wedge \cdots \wedge g^k(x)$.
Note that $\top (y(k)-x) \leq 0$. We claim that $\top (y(k)-x)= 0$.
Suppose not, so that $\top( y(k)-x) < 0$ for some $k$, 
which we may assume satisfies $k \geq 1$. It then follows that 
\begin{equation}
\label{eq-yk}
y(k) = g(x) \wedge \cdots \wedge g^k(x) \enspace. 
\end{equation}
Now choose $a > 0$ such that $\top(y(k) -x + ka) \leq 0$, 
which we may clearly do. Define $z(k) \in \R^n$ so that
\[ z(k) = x \wedge (g(x) + a) \wedge \cdots \wedge (g^{k-1}(x) + (k-1)a) \enspace. \]
Using (\ref{homo2}), (\ref{mono2}), we see that
\begin{gather*}
 g(z(k)) \leq g(x) \wedge \cdots \wedge (g^k(x) + (k-1)a)\enspace .
\end{gather*}
Since $a>0$, the right hand side is dominated
by $y(k)+(k-1)a$. Hence,
\begin{align*}
 g(z(k)) &\leq
\big( g(x) \wedge \cdots \wedge (g^k(x) + (k-1)a)\big)\wedge (y(k)+(k-1) a)\\
&= g(x) \wedge \cdots \wedge (g^{k-1}(x) + (k-2)a) \wedge (y(k) + (k-1)a) 
\enspace. \end{align*}
By choice of $a$, $y(k) + (k-1)a \leq x-a$. It follows that
\[ g(z(k)) \leq (x-a) \wedge g(x) \wedge \cdots \wedge (g^{k-1}(x) + (k-2)a) = z(k) - a \enspace. \]
But now, Proposition~\ref{prop-collatz-wielandt} implies that $\overline{\mychi}(g) \leq \top\big(g(z(k)) - z(k)\big) \leq -a < 0$, which is a contradiction. Hence, 
$\top (y(k) -x)= 0$ for all $k \in \N$, as claimed.

Since $y(k)-x$ is sequence of nonpositive vectors, each of which has at least one component $0$, there must be at least one coordinate, $1 \leq i \leq n$, such that $(y(k) - x)_i = 0$ for infinitely many $k \in \N$. But then, since $y(k) - x$ is decreasing, $(y(k) - x)_i = 0$ for all $k \in \N$. Hence, $g^k_i(x) \geq y(k)_i = x_i$, from which (\ref{e-new}) follows.
\end{proof}

\begin{cor}[Conjecture~2.1 of \cite{ctnxp}]
\label{cor-gk}
Let $f:\R^n\to\R^n$ be a topical function. There exists $1\leq i\leq n$ such that 
\begin{eqnarray*} 
\lim_{k\to \infty} f_i^k(y)/k &=& \overline{\mychi}(f)
\enspace ,
\end{eqnarray*}
for all $y\in \R^n$.
\end{cor}

\begin{proof} 
Since $f^k_i(x)/k \leq \top f^k(x)/k$, this follows from Theorem~\ref{theo-dyn} by letting $k \ra \infty$ and using (\ref{eq-scales}).
\end{proof}

\section{Existence of eigenvectors}

\subsection{Eigenvectors and bounded orbits}
\label{s2-evbo}
To prove Theorem~\ref{theo-orbits}, we need the result of an earlier paper, \cite[Lemma~4.2]{mxpnxp0}. We give the proof, for completeness. The result is a variant of Theorem~\ref{theo-orbits} for the supremum norm, where we take advantage of the lattice structure of $\R^n$. We note, following the comments before Theorem~\ref{theo-orbits}, that this is similarly a special case of a general result of Nussbaum in the context of nonexpansive functions on cones, \cite[Theorem~4.3]{nuss88}. 
\begin{lemma}
\label{lem-prev}
Let $f: \R^n \ra \R^n$ be a topical function. $f$ has an eigenvector with eigenvalue $h$ if, and only, if there exists $x \in \R^n$ such that $f^k(x) = kh + O(1)$ as $k \ra \infty$. (That is, $\supnorm{f^k(x) - kh}$ is bounded as $k \ra \infty$.)
\end{lemma}
\begin{proof}
If $f(x) = x + h$, then by (\ref{homo2}), $f^k(x) - kh = x$ for all $k \in \N$, which shows that the conclusion is necessary. Now suppose that the conclusion is satisfied and let $g = f - h$, so that $g^k(x)$ is bounded in the supremum norm as $k\to \infty$. Let
\[
u=\lim_{k\to\infty} \bigwedge_{\ell \geq k} g^{\ell}(x) 
\]
where the finiteness of $u$ follows from the boundedness of $g^k(x)$. By continuity and monotonicity of $g$,
\[
g(u)=\lim_{k\to\infty}
g(\bigwedge_{\ell\geq k} g^{\ell}(x))
\leq \lim_{k\to\infty}
\bigwedge_{\ell\geq k} g^{\ell+1}(x)
= u \enspace. 
\]
It follows from (\ref{mono2}) that $g^k(u)$ is nonincreasing as $k \ra \infty$. Since $g^k(x)$ is bounded as $k\to\infty$, it follows from (\ref{eq-scales}), that $g^k(u)$ is bounded too, so that $g^k(u)$ must converge to a limit, $v \in \R^n$. But then, by continuity of $g$, $g(v) = v$, so that $f(v) = h + v$, as required.
\end{proof}

Lemma~\ref{lem-prev} leads to the additive version of Theorem~\ref{theo-orbits}.

\begin{theorem}
\label{th-orbits}
Let $f: \R^n \ra \R^n$ be a topical function. $f$ has an eigenvector in $\R^n$ if, and only if, some (and hence all) orbits of $f$ are bounded in the Hilbert semi-norm.
\end{theorem}
\begin{proof}
Suppose that the orbit $\{ f^k(0) \;|\; k \in \N \}$ is bounded in the Hilbert semi-norm, so that $\hilbert{f^k(0)} \leq M$ for all $k \in \N$, for some $M > 0$. We show that $f$ has an eigenvector in $\R^n$. As discussed in \S\ref{s2-hilbert}, the rest of the argument is clear. Let $g = f - \overline{\mychi}(f)$, so that $\overline{\mychi}(g) = 0$. It follows from the definition of the Hilbert semi-norm that $\hilbert{g^k(0)} = \hilbert{f^k(0)} \leq M$ for all $k \in \N$. By (\ref{eq-mychik}), $\overline{\mychi}(g^k) = 0$ and so $\underline{\mychi}(g^k) \leq 0$. Applying Proposition~\ref{prop-collatz-wielandt} to $g^k$ and considering $x = 0$ in (\ref{eq-collatz}), we see that $\top g^k(0) \geq \overline{\mychi}(g^k)$. By symmetry, $\bot g^k(0) \leq \underline{\mychi}(g^k)$, so that $\top g^k(0) \geq 0 \geq \bot g^k(0)$. Note that if $a,b \geq 0$, then $a \vee b \leq a + b$. It follows that
\[ \supnorm{g^k(0)} = \top g^k(0) \vee -\bot g^k(0) \leq \top g^k(0) - \bot g^k(0) = \hilbert{g^k(0)} \leq M \enspace. \]
Hence $g^k(0)$ satisfies the conditions for Lemma~\ref{lem-prev} with $h = 0$ and so $g(x) = x$ for some $x \in \R^n$. It follows that $f(x) = x + \overline{\mychi}(f)$, as required. 
\end{proof}

We note that this proof does not use the full force of the Generalised Collatz-Wielandt formula (Proposition~\ref{prop-collatz-wielandt}); it requires Lemma~\ref{lemma-Lambda} but not Lemma~\ref{lemma-byk}.

\subsection{The associated graph}
\label{s2-graph}

We adapt the notation of (\ref{eq-uj}) to suit the additive framework. If $J \subseteq \{1, \cdots, n\}$, let $e_J \in \R^n$ denote the characteristic vector of $J$:
\begin{equation}
\label{eq-ujp}
(e_J)_i = \left\{\begin{array}{ll}
                 1 & \mrm{if $i \in J$} \\
                 0 & \mrm{if $i \not\in J$} \enspace.
                 \end{array}\right.
\end{equation}
Let $f:\R^n\to\R^n$ be a topical function. Define the associated graph of $f$, $\G(f)$, to be the directed graph with vertices $1, \cdots, n$ and an edge from $i$ to $j$, which we denote $i \ra j$, if, and only if, $\lim_{u \ra \infty}f_i(u e_{\{j\}})=\infty$. (The reader may care to check that if $f: (\Rp)^n \ra (\Rp)^n$ is a homogeneous, monotone function then $\G(f)$, as defined in \S\ref{s2-pft}, and $\G(\E(f))$, as defined here, are identical.) What follows is the additive version of Theorem~\ref{theo-main-new3} of \S\ref{s2-hilbert}.

\begin{theorem}
\label{th-new3}
If $f:\R^n \to \R^n$ is a topical function such that $\G(f)$ is strongly connected, then all super-eigenspaces of $f$ are bounded in the Hilbert semi-norm.
\end{theorem}
The proof of this relies on the following construction. For each edge $i\to j$ of $\G(f)$, define $h_{ji}: \Rim \to \R\cup\{-\infty\}$ by
\[
h_{ji}(x) =\sup\set{u \in \R}{f_i(u e_{\{j\}}) \leq x} \enspace,
\]
with the usual convention that $\sup \emptyset = -\infty$. (So that $h_{ji}(-\infty) = -\infty$.) Note that $h_{ji}(x) < \infty$ because there is an edge $i \ra j$. For any $\lambda\in \R$, we set $h^{\lambda}_{ji}(x)=h_{ji}(\lambda+x)$. 

\begin{lemma}
\label{lem-new}
Let $f:\R^n\to\R^n$ be a topical function and let $\lambda \in \R$. Let $i=i_1 \ra i_2 \ra \cdots \ra i_k=j$ be any directed path from $i$ to $j$ in $\G(f)$. Then, for all $x \in \R^n$
\begin{equation}
\label{e-bound}
(f(x)\leq \lambda+x\;\mrm{and}\; x\geq 0) \implies x_j \leq h^\lambda_{i_k i_{k-1}} \comp \cdots \comp h^\lambda_{i_2 i_1}(x_i) \enspace. 
\end{equation}
\end{lemma}
\begin{proof} Choose $x \in \R^n$ and $\lambda \in \R$ satisfying the conditions of (\ref{e-bound}). Since $x \geq 0$, it follows from (\ref{eq-ujp}) that $x_{i_p}e_{\{i_p\}} \leq x$, for all $1\leq p\leq k$. Hence, 
\[ f_{i_{p-1}}(x_{i_p}e_{\{i_p\}}) \leq f_{i_{p-1}}(x) \leq \lambda + x_{i_{p-1}}
\]
and so $x_{i_p} \leq h^\lambda_{i_p i_{p-1}}(x_{i_{p-1}})$. Putting these together, we deduce (\ref{e-bound}).
\end{proof}

\begin{proof}[Proof of Theorem~\ref{th-new3}] Since $\G(f)$ is strongly connected, we may choose, for each pair of vertices $i \not= j$, a directed path, $i=i_1\to\cdots \to i_k=j$, from $i$ to $j$ in $\G(f)$. Now choose $\lambda$ so that $S^\lambda(f)$ is nonempty and choose $x \in S^\lambda(f)$. Let $y = x - \bot x$. Note that $y \in S^\lambda(f)$, $\hilbert{y}=\hilbert{x}$ and $\bot y= 0$. It follows from (\ref{def-suphil}) that $\hilbert{y} = \supnorm{y}$. Since $\bot y = 0$, we may choose $i$ such that $y_i=0$. Let $j$ be any other coordinate. By (\ref{e-bound}), using the chosen path from $i$ to $j$, $0 \leq y_j\leq h^\lambda_{i_k i_{k-1}} \comp \cdots \comp h^\lambda_{i_2 i_1}(0)$. It follows that $\supnorm{y}$ is bounded for all $y$. Hence, $\hilbert{x} = \hilbert{y} = \supnorm{y}$ is also bounded, for all $x \in S^\lambda(f)$, from which the result follows.
\end{proof}

As mentioned in the Introduction, Amghibech and Dellacherie have introduced a graph, ${\mathcal Gr}(f)$, along the same lines as $\G(f)$. However, their construction uses a two-sided limit condition for an edge from $i$ to $j$:
\[ \lim_{u \ra \infty} f_i(ue_{\{j\}}) = \infty \;\;\mbox{and}\; \lim_{u \ra -\infty} f_i(ue_{\{j\}}) = -\infty \enspace.\]
If ${\mathcal Gr}(f)$ is strongly connected, so that $\G(f)$ is also strongly connected, Amghibech and Dellacherie make several deductions about the modified eigenvalue problem, $f(x) = u + x + \lambda$, where $u \in \R^n$ is an arbitrary constant. By using the fact that $f(x)-u$ is also a topical function, their conclusions follow from Theorem~\ref{theo-new} and Proposition~\ref{prop-collatz-wielandt}.

If $f$ is convex, the graph $\G(f)$ may be constructed more simply. Recall that a function $h: \R^n\to\R$ is convex if, for all $x,y \in \R^n$, $h(\lambda x + \mu y) \leq \lambda h(x) + \mu h(y)$, where $0 \leq \lambda,\mu \leq 1$ and $\lambda+\mu = 1$. A function $f: \R^n\to\R^n$ is convex if each component function $f_i: \R^n\to\R$ is convex. A simple deduction, which is left to the reader, captures the intuition that the derivative of $h$ is increasing. With the same notation, let $x' = \lambda x + \mu y = x+\mu(y-x) = y-\lambda(y-x)$. Then, 
\begin{equation}
\frac{h(x') - h(x)}{\mu} \leq \frac{h(y) - h(x')}{\lambda} \enspace.
\label{con2}
\end{equation}
For any function $f:\R^n\to\R^n$ define its {\em syntactic} graph, $G_s(f)$, to be the directed graph with vertices $1, \cdots, n$ and an edge from $i$ to $j$ if, and only if, $f_i$ depends on $x_j$ in the following sense: there is no map $h:\R^{n-1}\to\R$ such that $f_i(x)=h(x_1,\ldots,x_{j-1},x_{j+1},\ldots,x_n)$.

\begin{proposition}
\label{prop-convex}
Let $f:\R^n\to\R^n$ be a convex topical function. Then $\G(f)$ is identical to $\G_s(f)$. 
\end{proposition}
\begin{proof}
$\G(f)$ and $\G_s(f)$ have the same vertices, so it suffices to show that they have the same edges. It is clear that an edge of $\G(f)$ is also an edge of $\G_s(f)$. Conversely, suppose there is an edge from $i$ to $j$ in $\G_s(f)$. Then we can find $x,x'\in \R^n$ such that $x_k=x'_k$ for all $k\neq j$, $x_j\neq x'_j$, and $f_i(x) \neq f_i(x')$. Without loss of generality, assume that $x'_j>x_j$, so that $f_i(x') > f_i(x)$. Let $u > 0$ and let $y = x' + u e_j$. We may find $\lambda, \mu$, satisfying the convexity conditions, such that $x' = \lambda x + \mu y$. Rewriting (\ref{con2}), we see that
\[ f_i(x'+ ue_j) \geq \frac{u}{x'_j-x_j} (f_i(x')- f_i(x)) + f_i(x') \enspace,
\]
so that $\lim_{u \ra \infty} f_i(x'+ue_j) = \infty$. Now $x'+ue_j \leq \top x' + ue_j$, so using (\ref{homo2}) and (\ref{mono2}) we see that $f_i(ue_j) \geq f_i(x'+ue_j) - \top x'$. It follows that $\lim_{u \ra \infty} f_i(ue_j) = \infty$ so that there is an edge from $i$ to $j$ in $\G(f)$.
\end{proof}

Proposition~\ref{prop-convex} applies in particular to a nonnegative matrix, $A$, acting on the positive cone, as in the Introduction. Since the function $\R^2\to\R$ which takes $x\mapsto\log(\exp(x_1)+\exp(x_2))$ is convex, it is not difficult to show that $\E(A)$ is a convex topical function.

\subsection{Indecomposability}
\label{s2-indecomp}
We next prove the additive version
of Theorem~\ref{theo-main-new4bis}.
Recall that $f$ is {\em decomposable} 
if there is a partition $I \cup J = \{1, \cdots, n\}$, $I \cap J = \emptyset$, such that, $\forall i \in I$,
\begin{equation}
\label{eq-dec}
\lim_{u \ra \infty} f_i(u e_J) < \infty \enspace,
\end{equation}
and that $f$ is {\em indecomposable} if it is not decomposable.
\begin{theorem}
\label{theo-indeco}
Let $f:\R^n\to\R^n$ be a topical function. All super-eigenspaces of $f$ are bounded in the Hilbert semi-norm if, and only if, $f$ is indecomposable.
\end{theorem}
\begin{proof}
If $f:\R^n\to\R^n$ is a topical function for which, for some $\lambda\in \R$, $S^{\lambda}(f)$ is unbounded in the Hilbert semi-norm, we can find a sequence $x(k) \in S^{\lambda}(f)$ such that, $\bot x(k)=0$ for all $k \in \N$ and $\lim_{k\ra\infty} \top x(k) = \infty$. Since $[0,\infty]^n$ is compact for the usual topology, we may, possibly after replacing $x(k)$ by a subsequence, assume that $x(k)$ converges in $[0,\infty]^n$. Let $I=\set{1\leq i\leq n}{\lim_{k\ra\infty} x(k)_i<+\infty}$, and $J=\{1,\ldots,n\}\setminus I$. By construction, $J\neq \emptyset$, and, since $\bot x(k)=0$ for all $k$,  $I\neq \emptyset$. Moreover, for all $i\in I$, $\lim_{u \to\infty} f_i(ue_J) \leq \lim_{k \to \infty} f_i(x(k)) < +\infty$, which shows that $f$ is decomposable.

Conversely, let us assume that $f$ is decomposable, so that~\eqref{eq-dec} holds for some non-trivial partition $I\cup J=\{1,\ldots,n\}$. Possibly after a reordering of indices, we may assume that $I=\{1,\ldots, p\}$, $J=\{p+1,\ldots,n\}$. Let $q = n-p$ and let $y \in \R^p$, $z \in \R^q$. We may write $f(x)= (g(y,z), \;h(y,z))$, where $g: \R^p \times \R^q \ra \R^p$ and $h: \R^p \times \R^q \ra \R^q$. Let $w=\lim_{z\to \infty} g(0,z)$. The decomposability of $f$ implies that $w < \infty$. Now choose $\lambda \in \R$. If $w \leq \lambda$, then, for all $z\in \R^q$, $g(0,z) \leq 0 + \lambda$. Moreover, taking any $\nu \in \R$ with $\nu\geq 0$, and letting $z=(\nu,\ldots,\nu)\in \R^q$, we see that $h(0,z )\leq h(\nu,z)=\nu+h(0,0) \leq z + h(0,0)$. If we now choose $w \vee h(0,0) \leq \lambda$, then $(0,z)\in S^{\lambda}(f)$, for all $\nu\geq 0$. Hence, $S^\lambda(f)$ is unbounded in the Hilbert semi-norm.
\end{proof}

\subsection{Aggregated graphs}
\label{s2-aggregated}
In this section, we define the sequence of aggregated graphs used in Theorem~\ref{theo-main-new4} to characterize indecomposable topical functions.

For a set $X$, we denote by $\cP(X)$ the set of subsets of $X$. Given a positive integer $n$, we define inductively $\cP_1=\{1,\ldots,n\}$, and, for $k\geq 2$, $\cP_k=\cP(\cP_{k-1})$. We also define inductively the map $\sigma: \cup_{k\geq 1}\cP_k\to \cP_2$, by $\sigma(i)=\{i\}$ for all $i\in \cP_1$, and, for $X\in \cP_k$ with $k\geq 2$, $\sigma(X)= \bigcup_{Y\in X} \sigma(Y)$. For notational convenience, if $X\in \cup_{k\geq 1} \cP_k$, we use $e_X$ in place of $e_{\sigma(X)}$. For instance, for $n=4$, $X=\{\{1,2\},\{4\}\}\in \cP_3$, $\sigma(X)=\{1,2,4\}$, and $e_X=(1,1,0,1)$.

Recall that a strongly connected component of a directed graph is an equivalence classes of vertices under the relation of {\em communication}, \cite[Definition~2-3.7)]{bp94}: $i$ communicates with $j$ if, and only if, either $i=j$ or there are paths in both directions, from $i$ to $j$ and from $j$ to $i$. If $X,Y \subseteq \{1, \cdots, n\}$ are strongly connected components then $Y$ is {\em accessible} from $X$ if either $X = Y$ or there is a path from some vertex of $X$ to some vertex of $Y$. Accessibility is a partial order on the strongly connected components.

If $f:\R^n\to\R^n$ is a topical function, we define inductively the directed graphs $\G^k(f)$, for $k\geq 1$, as follows. The graph $\G^1(f)$ is by definition $\G(f)$. For $k\geq 2$, the vertices of $\G^k(f)$ are the strongly connected components of $\G^{k-1}(f)$, and there is an edge from the strongly connected component $I$ to the strongly connected component $J$ if, and only if,
\begin{equation}
\label{eq-edge}
\exists i \in \sigma(I) \;\;\mbox{such that}\; \lim_{u\to\infty} f_i(u e_{J})= \infty \enspace. 
\end{equation}
If there is at least one strongly connected component of $\G^k(f)$ not reduced to a single element, then, the number of vertices of $\G^{k+1}(f)$ is strictly less than the number of vertices of $\G^{k}(f)$. Otherwise, $\G^{k+1}(f)$ is isomorphic to $\G^k(f)$. We conclude that there is a least integer $N\leq n$ such that, for all $k\geq N$,  $\G^{k}(f)$ are isomorphic as directed graphs. We set $\G^{\infty}(f)=\G^{N}(f)$. All the strongly connected components of $\G^{\infty}(f)$ have only one element. 

To understand this construction, it may be helpful to consider the topical function $f:\R^4\to\R^4$ below. 
\begin{equation}
f(x) = \left(\begin{array}{c}
x_1 \vee (x_2 \wedge x_3 \wedge x_4)\\
x_3 \wedge x_4 \\
(x_1+x_2+x_4)/3 \\
x_3 \vee x_4 
\end{array}\right) \enspace.
\label{eq-example}
\end{equation}
The sequence of graphs $\G^1(f), \cdots, \G^4(f)=\G^{\infty}(f)$ is shown below.
\begin{center}
\begin{picture}(0,0)%
\includegraphics{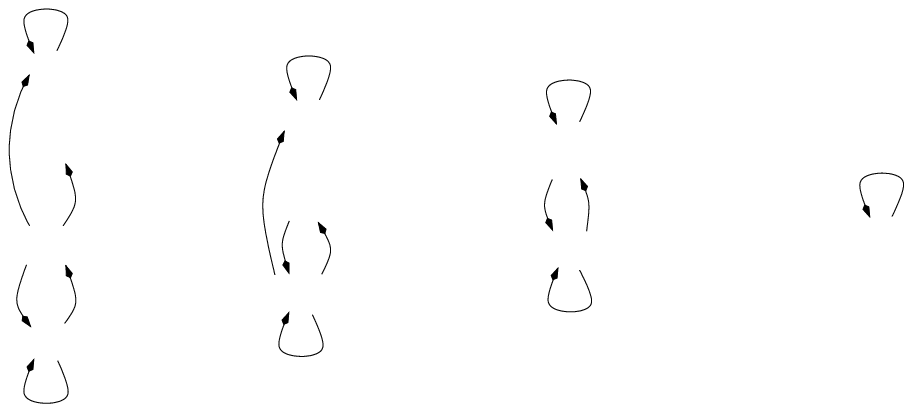}%
\end{picture}%
\setlength{\unitlength}{1973sp}%
\begingroup\makeatletter\ifx\SetFigFont\undefined%
\gdef\SetFigFont#1#2#3#4#5{%
  \reset@font\fontsize{#1}{#2pt}%
  \fontfamily{#3}\fontseries{#4}\fontshape{#5}%
  \selectfont}%
\fi\endgroup%
\begin{picture}(11544,3809)(458,-4628)
\put(8431,-3158){\makebox(0,0)[lb]{\smash{\SetFigFont{9}{10.8}{\rmdefault}{\mddefault}{\itdefault}{\color[rgb]{0,0,0}$\{\{\{1\}\},\{\{2\},\{3,4\}\}\}$}%
}}}
\put(458,-2798){\makebox(0,0)[lb]{\smash{\SetFigFont{9}{10.8}{\rmdefault}{\mddefault}{\itdefault}{\color[rgb]{0,0,0}$\G^1(f)$}%
}}}
\put(2596,-2776){\makebox(0,0)[lb]{\smash{\SetFigFont{9}{10.8}{\rmdefault}{\mddefault}{\itdefault}{\color[rgb]{0,0,0}$\G^2(f)$}%
}}}
\put(5461,-2776){\makebox(0,0)[lb]{\smash{\SetFigFont{9}{10.8}{\rmdefault}{\mddefault}{\itdefault}{\color[rgb]{0,0,0}$\G^3(f)$}%
}}}
\put(3859,-1936){\makebox(0,0)[lb]{\smash{\SetFigFont{9}{10.8}{\rmdefault}{\mddefault}{\itdefault}{\color[rgb]{0,0,0}$\{1\}$}%
}}}
\put(1456,-1471){\makebox(0,0)[lb]{\smash{\SetFigFont{9}{10.8}{\rmdefault}{\mddefault}{\itdefault}{\color[rgb]{0,0,0}$1$}%
}}}
\put(1441,-4111){\makebox(0,0)[lb]{\smash{\SetFigFont{9}{10.8}{\rmdefault}{\mddefault}{\itdefault}{\color[rgb]{0,0,0}$4$}%
}}}
\put(1456,-3196){\makebox(0,0)[lb]{\smash{\SetFigFont{9}{10.8}{\rmdefault}{\mddefault}{\itdefault}{\color[rgb]{0,0,0}$3$}%
}}}
\put(1456,-2311){\makebox(0,0)[lb]{\smash{\SetFigFont{9}{10.8}{\rmdefault}{\mddefault}{\itdefault}{\color[rgb]{0,0,0}$2$}%
}}}
\put(3814,-2716){\makebox(0,0)[lb]{\smash{\SetFigFont{9}{10.8}{\rmdefault}{\mddefault}{\itdefault}{\color[rgb]{0,0,0}$\{2\}$}%
}}}
\put(3769,-3631){\makebox(0,0)[lb]{\smash{\SetFigFont{9}{10.8}{\rmdefault}{\mddefault}{\itdefault}{\color[rgb]{0,0,0}$\{3,4\}$}%
}}}
\put(6158,-2258){\makebox(0,0)[lb]{\smash{\SetFigFont{9}{10.8}{\rmdefault}{\mddefault}{\itdefault}{\color[rgb]{0,0,0}$\{\{1\}\}$}%
}}}
\put(5888,-3218){\makebox(0,0)[lb]{\smash{\SetFigFont{9}{10.8}{\rmdefault}{\mddefault}{\itdefault}{\color[rgb]{0,0,0}$\{\{2\},\{3,4\}\}$}%
}}}
\put(7801,-2761){\makebox(0,0)[lb]{\smash{\SetFigFont{9}{10.8}{\rmdefault}{\mddefault}{\itdefault}{\color[rgb]{0,0,0}$\G^4(f)$}%
}}}
\end{picture}

\end{center}

The additive version of Theorem~\ref{theo-main-new4} requires the following lemma.
\begin{lemma}
\label{lem-nn}
Let $f:\R^n\to\R^n$ be a topical function and let $\lambda \in \R$. There is a function $H^{k}: \Rim \to \Rim$, such that, for all strongly connected components $X$ of $\G^k(f)$, for all $x \in \R^n$ and for all $i,j\in \sigma(X)$,
\begin{equation}
(f(x)\leq \lambda+x\;\mrm{and}\; x\geq 0) \implies x_j \leq H^{k}(x_i) \enspace .
\label{e-bound2}
\end{equation}
\end{lemma}
\begin{proof}
For $k = 1$ this follows from Lemma~\ref{lem-new}. There are only finitely many triples $(X,i,j)$, where $X$ is a strongly connected component of $\G(f)$ and $i,j \in X$. For each triple, either $i=j$ or we can find a path from $i$ to $j$ in $\G(f)$. In the latter case, Lemma~\ref{lem-new} gives a function $\Rim \ra \Rim$; in the former case, simply take the identity function. The pointwise maximum of these functions, over all triples, gives a function $H$ which satisfies (\ref{e-bound2}).

Now assume that $k\geq 2$. As before, there are finitely many triples of the form $(X,i,j)$, where $X$ is a strongly connected component of $\G^k(f)$ and $i,j \in \sigma(X)$. It is sufficient, by the same process of maximization as before, to choose a function $\Rim \ra \Rim$ satisfying (\ref{e-bound2}) for each triple $(X,i,j)$. Let $I$ and $J$ be the vertices of $\G^k(f)$ such that $i \in \sigma(I)$ and $j \in \sigma(J)$. If $I = J$, which includes the case $i=j$, then choose the function, $H = H^{k-1}$, provided by the inductive hypothesis.

Now assume that $I \not= J$. Then there is a path $I_1\to I_2 \to\cdots\to I_\ell$ in $\G^k(f)$, such that $I = I_1$ and $J = I_\ell$. By (\ref{eq-edge}) we can find $r_1 \in \sigma(I_1), \cdots, r_{\ell-1}\in \sigma(I_{\ell-1})$ such that, for $1 \leq m \leq \ell-1$,
\begin{equation}
\lim_{u \to \infty} f_{r_m} (u e_{I_{m+1}})= \infty \enspace .
\label{eq-arc}
\end{equation}
For all $\lambda \in \R$, $p \in\{1, \dots, n\}$ and $K\in \bigcup_{k\geq 1} \cP_k$, define the maps $h^{\lambda}_{K,p}:\Rim \to \R\cup\{\pm\infty\}$ by
\[
h^{\lambda}_{K,p}(t) = \sup\set{u \in \R}{f_p(u e_K) \leq \lambda + t} \enspace,
\]
with the usual convention that $\sup\emptyset=-\infty$. (So that $h^\lambda_{K,p}(-\infty) = -\infty$.) For $1 \leq m \leq \ell-1$, if $p=r_m$ and $K=I_{m+1}$, then by~\eqref{eq-arc}, $h^{\lambda}_{K,p}(t) < \infty$ for all $t$. We claim that if $f(x)\leq \lambda +x$ and $x\geq 0$ then 
\begin{equation}
x_j \leq H^{k-1} \comp h^{\lambda}_{I_\ell, r_{\ell-1}} \comp H^{k-1} \comp \cdots \comp H^{k-1} \comp h^{\lambda}_{I_2, r_1} \comp H^{k-1} (x_i) \enspace,
\label{eq-path2}
\end{equation}
which provides the required function $\Rim \ra \Rim$ for this triple. To conclude the proof, it suffices to show~\eqref{eq-path2}.

Choose the vertex $s_m\in \sigma(I_m)$ for $2 \leq m \leq \ell$, so that $x_{s_m}= \min_{v\in \sigma(I_m)} x_v$. Since $0\leq x$, we have $ x_{s_m}e_{I_m}\leq x$ and so
\[
f_{r_{m-1}}(x_{s_m}e_{I_m}) \leq f_{r_{m-1}}(x) \leq  \lambda + x_{r_{m-1}}\enspace.
\]
It follows that, for $2 \leq m \leq \ell$,
\begin{equation}
x_{s_m} \leq h^{\lambda}_{I_m,r_{m-1}} (x_{r_{m-1}}) 
\enspace .
\label{ineq-1}
\end{equation}
Moreover, by the inductive hypothesis, we have, for $2 \leq m \leq \ell-1$,
\begin{equation}
x_{r_m} \leq H^{k-1}(x_{s_m})
\label{ineq-2}
\end{equation}
while at the ends of the path, 
\begin{equation}
\label{ineq-3}
x_{r_1} \leq H^{k-1}(x_i) \;\;\mrm{and}\;\; x_j\leq H^{k-1}(x_{s_{\ell}}) \enspace.
\end{equation}
Composing the inequalities in~\eqref{ineq-1},~\eqref{ineq-2} and (\ref{ineq-3}) we get~\eqref{eq-path2}, as claimed.
\end{proof}

\begin{theorem}\label
{theo-main-new4b}
A topical function $f: \R^n \to \R^n$ is indecomposable if, and only if, $\G^{\infty}(f)$ is strongly connected. 
\end{theorem}

\begin{proof} Assume that $\G^\infty(f)$ is strongly connected and let $k$ be such that $\G^k(f)=\G^{\infty}(f)$. If $x \in S^\lambda(f)$ and $y = x - \bot x$ then, arguing as in the proof of Theorem~\ref{th-new3}, we can use Lemma~\ref{lem-nn} to show that $\hilbert{x} = \supnorm{y} \leq H^k(0)$. Hence, all nonempty super-eigenspaces of $f$ are bounded in the Hilbert semi-norm and so, by Theorem~\ref{theo-indeco}, $f$ is indecomposable.

Conversely, let us assume that $\G^\infty(f)= \G^k(f)$ is not strongly connected. Recall the partial order of accessibility defined above. Choose a vertex $X$ of $\G^k(f)$ which is minimal in this partial order. Let $J=\sigma(X)$ and $I=\{1,\ldots,n\} \setminus J$. For all $i\in I$, there is a strongly connected component $Y$ of $\G^k(f)$ such that $i \in \sigma(Y)$. Note that $Y \not= X$. If $\lim_{u\to\infty} f_i (u e_J) = \infty$, then by (\ref{eq-edge}) there would be an edge from $Y$ to $X$ in $\G^k(f)$, which contradicts the minimality of $X$. Therefore, for all $i \in I$, $\lim_{u\to\infty} f_i(u e_J) < \infty$, which implies by (\ref{eq-dec}) that $f$ is decomposable.
\end{proof}

As might be expected in the light of Proposition~\ref{prop-convex}, the aggregation process simplifies when $f$ is convex.

\begin{proposition}
\label{prop-convex2}
If $f:\R^n\to\R^n$ is a convex topical function then $\G^\infty(f)$ is isomorphic to $\G^2(f)$.
\end{proposition}
\begin{proof}
It suffices to show that the aggregation process stabilises by $k = 2$. Let $I$ and $J$ be strongly connected components of $\G(f)$ for which there is an edge $I \ra J$ in $\G^2(f)$. By (\ref{eq-edge}), $\exists i \in \sigma(I)$ such that $\lim_{u \ra \infty}f_i(ue_J) = \infty$. It follows that $\exists j \in \sigma(J)$ such that there is an edge $i \ra j$ in the syntactic graph, $\G_s(f)$, for otherwise, $f_i$ would not depend on any of the components $x_j$ for $j \in J$. By Proposition~\ref{prop-convex}, there is an edge $i \ra j$ in $\G(f)$. Since $I$ and $J$ are strongly connected components of $\G(f)$, there is then a path in $\G(f)$ from any vertex of $I$ to any vertex of $J$. It follows that $\G^2(f)$ cannot have any strongly connected components other than single vertices, so that the aggregation process must have stabilised by $k=2$. 
\end{proof}

\subsection{Slice spaces and recession functions}
\label{s2-slice2}
Let $f: \R^n \ra \R^n$ be a topical function. In the additive context, the {\em recession function} of $f$, $\hat{f}: \R^n \ra \R^n$ is defined by $\ha f(x)= \lim_{t\to\infty} t^{-1}f(t x)$. The recession function does not exist in general. For example, any function of the form $f:\R^2\to\R^2$, $(x_1,x_2)\mapsto (x_1, x_1+h(x_2-x_1))$ is topical provided $h$ has a derivative which satisfies $0\leq h'(t)\leq 1$, for all $t$. It is not difficult to find such an $h$ such that $\lim_{t \ra \infty} t^{-1}h(t)$ does not exist, so that $\lim_{t \ra \infty} t^{-1}f(tx)$ does not exist at $x = (0,1)$. Nevertheless, the recession function does exist for many interesting examples of topical functions, including all convex topical functions and all the other examples discussed in this paper.

It is not difficult to see that $\ha f$ is also a topical function and that it is multiplicatively homogeneous: $\ha f(t x)=t \ha f(x)$, $\forall t>0$. It follows that $\ha f(0) = 0$ so that $\ha f$ always has the {\em trivial eigenvectors}, $(u, \cdots, u)$, for $u \in \R$. The following is the additive version of Theorem~\ref{theo-main-new5}.

\begin{theorem}
\label{theo-fhat}
If a topical function $f$ has a recession function whose only eigenvectors are trivial, then all slice spaces of $f$ are bounded in the Hilbert semi-norm.
\end{theorem}
\begin{proof}
If $S^\lambda_\mu(f)$ is unbounded in the Hilbert semi-norm, then we can find a sequence $x(k)\in S^\lambda_\mu(f)$ such that $\bot x(k)=0$ for all $k \in \N$ and $\lim_{k \ra \infty} \top x(k)=+\infty$. We may assume that $\top x(k) > 0$ for all $k \in \N$. Since $y(k)=x(k)/\top x(k) \in [0,1]^n$, then, possibly after replacing $x(k)$ by a subsequence, we may assume that $y(k)\to y$ as $k \ra \infty$, for some $y\in [0,1]^n$. Note that $\bot y=0$ and $\top y=1$, so that $y$ is not of the form $(u, \cdots, u)$. By (\ref{nxp-supremum}),
\[ \supnorm{f(\top x(k).y(k)) - f(\top x(k).y)} \leq \top x(k) \supnorm{y(k) - y} \enspace. \]
Dividing by $\top x(k)$ and letting $k \ra \infty$, so that $\top x(k) \ra \infty$, we see that
\[ \lim_{k \ra \infty} (\top x(k))^{-1}f(\top x(k).y(k)) = \hat{f}(y)\enspace. \]
Since $x(k) \in S^\lambda_\mu(f)$, $\mu \leq f(x(k)) - x(k) \leq \lambda$. Dividing by $\top x(k)$ and letting $k \ra \infty$, we see that $\hat{f}(y) = y$. Hence $\hat{f}$ has a non-trivial eigenvector. The result follows.
\end{proof}

If $\ha f$ is a multiplicatively homogeneous topical function, it is easy to see that $\ha f$ has only trivial eigenvectors if, and only if, its eigenspace is bounded in the Hilbert semi-norm.

The condition of Theorem~\ref{theo-fhat} is not necessary for all the slice spaces of $f$ to be bounded. For example, let $h: \R\to \R$ be defined by $h(t)=1+t- \sqrt{1+t}$ and $h(-t) = -h(t)$ for $t\geq 0$. It is clear that $h(t)$ is differentiable and that $0 \leq h'(t) \leq 1$ for all $t \in \R$, so that the function $f:\R^2\to \R^2$, $(x,y) \mapsto (x,x+h(y-x))$, is topical. If $(x,y)\in S^{\lambda}_\mu(f)$, then $\mu+y \leq x+h(y-x) \leq \lambda +y$. It follows easily that $y-x$ is bounded, so that all slice spaces are bounded in the Hilbert semi-norm. However, $\hat f(x) = x$, so that every vector is an eigenvector.

\bibliographystyle{plain}

\end{document}